\documentclass[11pt]{amsart}
\usepackage{amsbsy,amsfonts}
\usepackage{latexsym,eucal,amsmath,amsthm,epsfig}
\usepackage{amssymb}
\theoremstyle{plain}
     \newtheorem{theorem}[subsection]{Theorem}

\theoremstyle{remark}

\theoremstyle{definition}

\begin{document}
\date{}
\title[Finite difference approximations for hyperbolic PDE with point-wise delay]{Finite difference approximations for the first-order hyperbolic partial differential
 equation with point-wise delay}
\author{Paramjeet Singh}
\address{Laboratoire Jacques-Louis Lions, Universit\'{e} Pierre et Marie Curie,
Paris, France}
\address{Department of Mathematics, Panjab University,
Chandigarh, India}

\author{Kapil K. Sharma}
\address{Department of Mathematics, Panjab University, Chandigarh, India}

\begin{abstract}
Explicit numerical methods based on Lax-Friedrichs and Leap-Frog finite difference
approximations are constructed to find the numerical solution of the first-order
hyperbolic partial differential  equation  with point-wise delay or advance, $i.e.$,
shift in  space.  The differential equation involving point-wise delay and advance
models  the distribution of the time intervals between successive neuronal firings.
In this paper, we continue the numerical study which was initiated in \cite{Singh2}.
We construct higher order numerical approximations and discuss their  consistency,
stability and convergence. The numerical approximations constructed in this paper
are consistent, stable under CFL condition, and convergent. We also extend our methods
to the higher space dimensions.  Some test examples are included to illustrate our
approach. These examples verify the theoretical estimates and shows the effect of
point-wise delay on the solution.
\end{abstract}
\maketitle

\section{Introduction}
Partial differential difference equations or more generally partial functional differential
equations  can be found in several mathematical models of control theory, mathematical
biology, climate models, mathematical economics, the theory of systems which communicate
through lossless channels, meteorology, and many other areas, see
\cite{Hale,Perthamebook,Wu,Kolmanovskii,Stein1965}.
Several biological phenomena can be mathematically  modeled by using the time dependent
first-order partial differential equations of hyperbolic type which contains point-wise
delay or shift in space. A detailed mathematical analysis and  numerical computations of
the ordinary delay differential equations is discussed by Bellen and Zennaro \cite{Bellen}
and references therein.
In this paper, we continue the study which was initiated in~\cite{Singh2}, $i.e.$,
we consider the following first-order hyperbolic partial differential  equation having
 point-wise delay with an initial data $u_0$ on domain $\Omega:=(0,X)$. In general it reads
\begin{equation}\label{1.1}
\begin{array}{rcl}
u_t+a u_x&=&b u(x-\alpha,t),
\quad  \mbox{$x\in \Omega$, $t>0$}\\
u(x,0) & = & u_0(x), \quad  x\in \bar\Omega\\
u(s,t)& = & \phi(s,t), \quad\ \mbox{$\forall \ s\in[-\alpha ,0]$, $t>0$},\quad \mbox {for $a>0$}\\
u(X,t)&= & \psi(t), \quad \mbox {for $a<0$},
\end{array}
\end{equation}
where $a=a(x,t)$ and $b=b(x,t)$  are sufficiently smooth functions of $x$ and $t$, $\alpha$
is the value of the point-wise delay which is non-zero fixed real number.  We also assume
that $a$ and $b$ are bounded in the entire domain and $a$ is not changing its sign in the
entire domain. Let $|a(x,t)|<A$ and  $|b(x,t)|<B$, $\forall (x,t)$. The unknown function
$u$ is  defined in the underlying domain and also in the intervals $[-\alpha,0]$ due to
the presence of point-wise delay.  So our  domain is $[-\alpha,0]\cup[0,X]$ and $t>0$.
The coefficients are sufficiently smooth functions in these intervals and the unknown
function $u$ is as smooth as the initial data. Due to the presence of point-wise delay
in equation (\ref{1.1}), we need a boundary-interval condition in the left side of domain,
$i.e.$, in the interval $[- \alpha ,0]$. The equation (\ref{1.1}) is first-order hyperbolic
 with difference terms, so it requires one boundary condition according to the direction
of characteristics~\cite{Morton}. A similar equation to (\ref{1.1}) having positive shift
in right side is obtained by replacing $-\alpha$ with $\beta$, see \cite{Singh2}.
Numerical schemes for equation  having positive shift can be obtained in similar
fashion by replacing point-wise delay with advance.  \\
Due to the presence  of point-wise delay   and non-constant coefficients, it is not
difficult but to some extent impossible  to find the analytic solution of such type of
partial differential equations by using the  usual methods to find the exact solution of
partial differential equations, see~\cite{Evans}. Also we cannot solve such type of partial
differential difference equations with classical numerical methods
\cite{Roos,Godlewski,Kulikovskii,Schroll}. Ordinary differential equations with
difference terms (delay) are quite well understood by now but there is no comparable
theory for partial differential equations (i.e. for time and space dependent unknowns).
Numerical solution of a general class of delay differential equation,  including stiff
problem, differential-algebraic delay equations and neutral problem is discussed by
Guglielmi et al. in \cite{Hairer1}. Implicit Runge-Kutta method is applied  in modified
form and possible difficulties are discussed. Asymptotic stabilities properties of
implicit Runge-Kutta method for ordinary delay differential equations  is considered by
Hairer et al. in \cite{Hairer2}. Parabolic partial differential equations with time delay
is considered in \cite{Gander}. The equation (\ref{1.1}) is of the form 1-D scalar
conservation laws with source term in right side which is the function of $x$ and $t$
and this source contains point-wise delay.\\
If delay and advance arguments are sufficiently small, the authors used the  Taylor
series approximations for the difference arguments and proposed an explicit numerical
scheme based on the upwind finite difference method which is discussed in~\cite{Singh1}.
This method has restriction on the size of point-wise delay and advance.  Equation containing
only point-wise delay or advance is considered by the authors in~\cite{Singh2}. This scheme
tackles small as well as large difference argument and has first-order convergence in both space and time. \\
In the following, we construct numerical approximations based on the Lax-Friedrichs and the
Leap-Frog finite difference approximations to increase the order of convergence. The applicability
of the second-order finite difference method remains valid to solve such type of hyperbolic partial
differential difference equations provided the geometry of underlying domain is not complicated.
Also it is easy to extend our ideas into higher space dimensions in finite difference approximations.
In the construction of numerical schemes,  a special type of mesh is generated during the discretization
so that  the difference argument also belong to the discrete set of grid points. The numerical
methods works  for large as well as small values of point-wise delay. We  construct the numerical
schemes to find the approximate solution of problem (\ref{1.1}) in section 2 and discuss the
consistency, stability and convergence.  Also we discuss the extensions of numerical methods
in higher space dimensions in section 3. In section 4,  we include some numerical test examples
to validate the predicted theory. Finally, in section 5, we make some concluding remarks
illustrating the effect of difference arguments on the solution behavior.

\section{Numerical Approximations}
In this section, we construct  numerical schemes based on the finite difference method
\cite{Morton,Leveque,Strikwerda}. We construct  explicit numerical approximations for
the given equation (\ref{1.1}) based on Lax-Friedrichs and Leap-frog finite difference
approximations. For space-time approximations based on finite differences, the $(x,t)$
plane is descretized by taking mesh width $\Delta x$ and  time step $\Delta t$, and
defining the grid points $(x_j,t_n)$ by
$$x_j= j \Delta x, \ j=0,1,...,J-1,J; \quad t_n =n \Delta t ,\ n=0,1,2....$$
Now we look for the approximate numerical solution $U_j^n$ that approximate $u(x_j,t_n), \forall j,n$.
We write the closure of $\Omega_{\Delta x}$ as $\bar\Omega_{\Delta x}$ and $\bar\Omega_{\Delta x}=(x_j=j \Delta x,\ j=0,1,2,...,J)$.
\subsection{Lax-Friedrichs approximation:}
In this numerical approximation, we approximate the time derivative by forward difference
and space by centered difference and then we replace $U_j^n$ by the mean value between
$U_{j+1}^n$ and $U_{j-1}^n$ for stability purpose, see \cite{Morton}.
Numerical approximation  is given by
\begin{equation}\label{2.1}
\frac{U_j^{n+1}-\frac{U_{j+1}^n+U_{j-1}^n}{2}}{\Delta t}+a_j^n \frac{U_{j+1}^n-U_{j-1}^n}{2\Delta x}=b_j^n u(x_j-\alpha,t_n).
\end{equation}
To tackle the point-wise delay  in the numerical approximation (\ref{2.1}),
we discretize the domain in such a way that  $(x_j-\alpha)$ is  a nodal point,
$\forall~j=0,1,...,J$, $i.e.$, we choose $\Delta x~s.t.~\alpha = m_0 \Delta x ,m_0 \in \mathbb{N}$
and we take total number of points in $x-$direction s.t.
$$ J=\frac{X}{\Delta x}=k X \frac{\mathrm{mantissa}(\alpha)}{\alpha}, ~k \in \mathbb{N},$$ where
mantissa of any real number is defined as  positive fractional part of that number.\\
The term containing point-wise delay $(\forall~j=0,1,...,J)$  can be written as\\
$\quad \quad \quad \quad \quad  u(x_j -\alpha,t_n)~=u(j\Delta x -m_0 \Delta x,t_n)$  \\
$\quad \quad \quad \quad\quad \quad \quad \quad  \quad\quad\quad             =u((j -m_0 )\Delta x,t_n)$ \\
 $\quad\quad\quad\quad  \quad\quad\quad\quad\quad\quad  \quad             \approx U_{j-m_0 }^n.
$\\
Therefore, the numerical approximation is given by
\begin{equation}\label{2.2}
\frac{U_j^{n+1}-\frac{U_{j+1}^n+U_{j-1}^n}{2}}{\Delta t}+a_j^n \frac{U_{j+1}^n-U_{j-1}^n}{2\Delta x}=b_j^n  U_{j-m_0}^n,  \quad \forall\ j=1,2,...,J-1,
\end{equation}
together with initial  and boundary-interval  conditions as following
\begin{subequations}\label{2.2a}
\begin{align}
&U_j^0=u^0(x_j),\quad  \forall \ j=1,2,...,J-1\\
&U_0^n=u(s,t_n)= \phi(s,t_n),~\forall \ s \in[- \alpha ,0], ~n=0,1,2\dots\\
&U_J^n= \psi(t_n), \ n=0,1,2\dots.
\end{align}
\end{subequations}
\textbf{Now we prove consistency, stability, and convergence of the
numerical approximation  (3-4).}
\begin{theorem}\label{theorem:2.1}
Let the function $u_0$ be sufficiently smooth then the explicit numerical scheme
(\ref{2.2}) is consistent of order 2 in space and of order 1 in time. Also if CFL
condition is satisfied, then on the  arbitrary finite time interval $[0,t_f]$, the
numerical scheme (\ref{2.2}) is stable in the maximum norm where stability constant
is of the form $C=1+O(\Delta t)$, and hence convergent.
\end{theorem}
\begin{proof}
Firstly, we study the consistency of the numerical scheme (\ref{2.1}). The consistency
error of the numerical  scheme (\ref{2.2}) is the difference between the two sides of the
equation when the approximation $U_j^n$ is replaced throughout by the exact solution $u(x_j,t_n)$
of the differential equation (\ref{1.1}). If $u$ is sufficiently smooth, then the consistency error
$T_j^n$ of the difference scheme (\ref{2.2}) is given by
$$T_j^n =\frac{u(x_j,t_{n+1})-\frac{u(x_{j+1},t_n)+u(x_{j-1},t_n)}{2}}{\Delta t}+ a_j^n\frac{u(x_{j+1},t_n)-u(x_{j-1},t_n)}{2\Delta x}-b_j^n u_{j-m_0}^n. $$
By using  the Taylor series approximation for the term $u(x_j,t_{n+1})$ w.r.t. to $t$ and
for the  terms $u(x_{j-1},t_{n})$ and $u(x_{j+1},t_{n})$   w.r.t. $x$, we get
\begin{align*}
T_j^n & = [u_t+\frac{\Delta t}{2} u_{tt}-\frac{\Delta x^2}{\Delta t}u_{xx}+O(\Delta t^2)+O(\frac{\Delta x^4}{\Delta t^2})]_j^n+[a(u_x+\frac{\Delta x^2}{6} u_{xxx})+O(\Delta x^4)]_j^n-b_j^n u_{j-m_0}^n\\
&=[u_t+au_x]_j^n-b_j^n u_{j-m_0}^n+\frac{\Delta t}{2} u_{tt}-\frac{\Delta x^2}{\Delta t}u_{xx}+\frac{\Delta x^2}{6}a_j^n u_{xxx}+O(\Delta x^4+\Delta t^{-1}\Delta x^4+\Delta t^2).
\end{align*}
As $u$ satisfy the given differential  equation (\ref{1.1}), we have $$\quad [u_t+au_x]_j^n-b_j^n u_{j-m_0}^n=0.$$
Therefore, the consistency error is given by
\begin{equation*}
T_j^n =\frac{\Delta t}{2} u_{tt}-\frac{\Delta x^2}{\Delta t}u_{xx}+\frac{\Delta x^2}{6}a_j^n u_{xxx}+O(\Delta x^4+\Delta t^{-1}\Delta x^4+\Delta t^2).
\end{equation*}
Therefore,  $T_j^n\rightarrow 0$ while $(\Delta x,\Delta t)\rightarrow(0,0)$, which implies the
numerical scheme is consistent of order 2 in space and of order 1 in time as long as $\Delta t^{-1}\Delta x^2\longrightarrow 0$.\\
Now by solving the   finite difference scheme (\ref{2.2}) for $U_j^{n+1}$, we get
\begin{equation}\label{2.4}
U_j^{n+1}=\frac{1}{2}(U_{j+1}^{n}+U_{j-1}^n)-\frac{a_j^n\Delta t}{2\Delta x}(U_{j+1}^n-U_{j-1}^n)+\Delta t b_j^n U_{j-m_0}^n, \quad \forall\ j=1,2,...,J-1.
\end{equation}
By applying the triangle inequality, we obtain
$$|U_j^{n+1}|\leq \frac{1}{2}\left|\left(1-\frac{a_j^n \Delta t}{\Delta x}\right)\right|~|U_{j+1}^n|+ \frac{1}{2}\left|\left(1+\frac{a_j^n \Delta t}{\Delta x}\right)\right|~|U_{j-1}^n|+\Delta t|b_j^n||U_{j-m_0}^{n}.$$
Now by taking the maximum norm, we get
\begin{align*}
\|U^{n+1}\|_{\infty,\Delta x}&= \max_j |U_j^{n+1}| \\
                    & \leq  \frac{1}{2}\max_j\left|\left(1-\frac{a_j^n \Delta t}{\Delta x}\right)\right|~|U_{j+1}^n|+ \frac{1}{2}\max_j\left|\left(1+\frac{a_j^n \Delta t}{\Delta x}\right)\right|~|U_{j-1}^n|+\Delta t\max_j|b_j^n|~|U_{j-m_0}^{n}|.
 \end{align*}
Using CFL condition $\frac{A \Delta t}{\Delta x} \leq1$ [where $A$ is the bound of $a(x,t)$, $\forall\ (x,t)]$, first two terms in above inequality can be combined and we get
\begin{align*}
 \|U^{n+1}\|_{\infty,\Delta x} & \leq (1+B \Delta t)\|U^n\|_{\infty,\Delta x},
\end{align*}
where $|b(x,t)|<B$, $\forall \ (x,t)$. The term $B \Delta t$  can be controlled by $\Delta t$ from which we can predict that  effect of the term  $B \Delta t$  is of the form  $O(\Delta t)$. Using these values, we get the following estimate
$$ \|U^{n+1}\|_{\infty,\Delta x} \leq  (1+O(\Delta t))\|U^n\|_{\infty,\Delta x}$$
$i.e.,$ $$ \|U^{n+1}\|_{\infty,\Delta x} \leq  C\|U^n\|_{\infty,\Delta x},$$
which implies the stability of the numerical approximation, where stability constant $C$ is of the form $C=1+O(\Delta t)$.\\
Now the error in the numerical approximation (\ref{2.4}) is given by
$$ e_j^n=U_j^n-u(x_j,t_n).$$
We set $e_j^n=U_j^n-u(x_j,t_n)$ in (\ref{2.4}).  The approximate solution $U_j^n$ satisfies the difference equation (\ref{2.4}) exactly, while exact solution $u(x_j,t_n)$ leaves the remainder $T_j^n \Delta t$.  Therefore, the error in the numerical approximation is given by
$$e_j^{n+1}=(1-a_j^n \frac{\Delta t}{\Delta x})e_{j+1}^n+(1+a_j^n \frac{\Delta t}{\Delta x}) e_{j-1}^n+b_j^n \Delta t~e_{j-m_0}^n-\Delta t T_j^n, $$
and  \quad \qquad  \quad $e_0^n=0$.\\
Let $ E^n = \max_j \{|e_j^n|,  j=0,1,..., J\}$.\\
For $|a_j^n \frac{\Delta t}{\Delta x}|\leq 1$,
\begin{align*}
 E^{n+1} &= \max_j |e_j^{n+1}|\\
         & \leq E^n+|b_j^n| \Delta t E^n+ \Delta t \max_j |T_j^n|\\
          & \leq E^n+B \Delta t E^n+ \Delta t \max_j |T_j^n| \\
          & =(1+B \Delta t)E^n +\Delta t \max_j |T_j^n|,
\end{align*}
since we are using the given initial value for $U_j^n$, so $E^0=0$ and if we suppose that the  consistency error is bounded $i.e.$ $|T_j^n|\leq T_{max},$ then by using induction method in the above inequality
$$E^n\leq n\Delta t T_{max}\leq t_f T_{max},$$ where $n\Delta t=t_f$,
which proves that the numerical scheme (\ref{2.2}) is convergent provided that the solution $u$ has bounded derivatives up to second order.
\qed
\end{proof}
%\quad\quad \quad\quad$\Box$ Q.E.D.\\

\subsection{Leap-Frog approximation:}
In this numerical approximation, we use central difference for both the space and time.  The numerical approximation for equation (\ref{1.1}) is given by
\begin{equation}\label{2.5}
\frac{U_j^{n+1}-U_{j}^{n-1}}{2\Delta t}+a_j^n\frac{U_{j+1}^n-U_{j-1}^n}{2\Delta x}=G_j^n, \quad \forall \ j=1,2,...,J-1;
\end{equation}
together with initial and boundary-interval  conditions as given in (\ref{2.2a}).\\
We write the numerical approximation as below
\begin{equation}\label{2.6}
U_j^{n+1}=U_{j}^{n-1}-\frac{a_j^n \Delta t}{\Delta x}(U_{j+1}^n-U_{j-1}^{n})+2\Delta t b_j^n U_{j-m_0}^{n}, \quad \forall \ j=1,2,...,J-1.
\end{equation}
This scheme uses two time intervals to get a central time difference and spreads its legs to pick up the space difference at the intermediate time level. It is an explicit scheme which requires a special technique to get it started. The initial condition will usually determine the values of $U^0$ but $U^1$ can be obtained by any convenient one-step scheme. In this case, we initialize the numerical scheme (\ref{2.6}) with the Lax-Friedrichs scheme (\ref{2.2}). Then this scheme gives $U^2$, $U^3$,... in succession.\\
 Now we discuss consistency, stability and convergence of this numerical approximation.  The consistency error of the numerical  scheme (\ref{2.5}) is the difference between the two sides of the equation when the approximation $U_j^n$ is replaced throughout by the exact solution $u(x_j,t_n)$ of the differential equation (\ref{1.1}). If $u$ is sufficiently smooth, then the consistency error $T_j^n$ of the difference scheme (\ref{2.5}) is given by
\begin{align*}
T_j^n &=\frac{u(x_j,t_{n+1})-u(x_j,t_{n-1})}{2\Delta t}+ a_j^n\frac{u(x_{j+1},t_n)-u(x_{j-1},t_n)}{2\Delta x}-b_j^n u_{j-m_0}^n\\
        &= [u_t+\frac{\Delta t^2}{6} u_{ttt}+O(\Delta t^4)]_j^n+[a(u_x+\frac{\Delta x^2}{6} u_{xxx})+O(\Delta x^4)]_j^n-b_j^n u_{j-m_0}^n\\
          &=[u_t+au_x]_j^n-b_j^n u_{j-m_0}^n +\frac{\Delta t^2}{6} u_{ttt}+\frac{\Delta x^2}{6} u_{xxx}+O(\Delta x^4)+O(\Delta t^4).
\end{align*}
Here we used the Taylor series approximations for the terms $u(x_j,t_{n+1})$ and $u(x_{j-1},t_{n})$ w.r.t. $t$ and $x$ respectively.\\
As $u$ satisfy the given differential  equation (\ref{1.1}),$$\quad [u_t+au_x]_j^n-b_j^n u_{j-m_0}^n=0.$$
Therefore, the consistency error is given by
\begin{align*}
T_j^n &= \frac{\Delta t^2}{6} u_{ttt}+\frac{\Delta x^2}{6} u_{xxx}+O(\Delta x^4)+O(\Delta t^4).
\end{align*}
Now $T_j^n\rightarrow 0$ while $(\Delta x,\Delta t)\rightarrow(0,0)$, which implies that the numerical scheme (\ref{2.5}) is consistent of order 2 in both space and and time.\\
Now the Leap-Frog scheme (\ref{2.5}) is consistent of order 2 in both space and time, the stability condition is the same as in the case of right hand side is zero, see \cite{Strikwerda}. When the term containing point-wise delay is zero, it has been shown \cite{Strikwerda} that Leap-Frog scheme is stable if the CFL number is strictly less than  1, that is $A\frac{\Delta t}{\Delta x}< 1$, where $A$ is the bound of $a$. Also we have proved in the  previous numerical scheme (\ref{2.2}) that right hand side effects only on the form of stability constant. Therefore, the numerical approximation is  stable in the $L_2-$norm provided  CFL condition is satisfied.  Now the proposed finite difference approximation  is linear, therefore,  by Lax-Richtmyer Equivalence Theorem it is convergent~\cite{Strikwerda}.
\section{Extension to Higher Spatial Dimensions}
We now consider the extensions of our numerical schemes to the higher spatial dimensions. For the sake of simplicity, we consider the problem in two spatial dimensions. The extension to three spatial dimensions can be done in a similar fashion \cite{Lax,Zhang}. The natural generalization of the one-dimensional  model problem  (\ref{1.1}) is the following equation together with initial data and boundary-intervals conditions
\begin{equation}\label{3.1}
u_t+a u_x+bu_y=c u(x-\alpha,y-\beta,t),
\end{equation}
where $a$, $b$, $c$ are functions of $x$, $y$ and $t$. $\alpha$ and  $\beta$ are the values of the point-wise delay in $x$ and $y$-direction respectively. We consider the rectangular domain in the $(x,y)-$plane as $0<x<X, \quad 0<y<Y$. For numerical approximations,
we discretize  the domain by taking uniform grid points with a spacing $\Delta x$ in the $x-$direction and $\Delta y$ in the $y-$direction. The grid points $(x_i,y_j,t_n)$ are defined as following
$$x_i= i \Delta x, \ i=0,1,...,J_x; \quad y_j= j \Delta y, \ j=0,1,...,J_y; \quad t_n =n \Delta t ,\ n=0,1,2....$$
Now we write the extension of Lax-Friedrichs scheme (\ref{2.2}). The approximate solution is denoted by  $U_{i,j}^n$. Numerical scheme is given by
\begin{equation}\label{3.2}
U_{i,j}^{n+1}=\frac{1}{4}(U_{i+1,j}^n+U_{i-1,j}^n+U_{i,j+1}^n+U_{i,j-1}^n)-\frac{\Delta t}{2\Delta x}a_{i,j}^n(U_{i+1,j}^n-U_{i-1,j}^n)
\end{equation}
$$-\frac{\Delta t}{2\Delta y}b_{i,j}^n(U_{i,j+1}^n-U_{i,j-1}^n)+\Delta t c_{i,j}^n U_{i-m_0,j-q_0}^n,$$
together with appropriate initial data and boundary-interval conditions. Here we take the grid points in both the directions $(x \ \mathrm{and} \ y)$ in such a way that the term containing point-wise delays is also belong to discrete set of  grid points which can be done as we did in the one dimensional  case. We take total number of  points in both $x$ and $y$ direction such that corresponding delays are on $m_0$ and $q_0$ node  and total number of points in both the directions  are given by
$$\quad J_x=\frac{X}{\Delta x}=k X \frac{\mathrm{mantissa}(\alpha)}{\alpha}, ~k \in \mathbb{N},$$
$$\mathrm{and} \ J_y=\frac{Y}{\Delta y}=r Y \frac{\mathrm{mantissa}(\beta)}{\beta}, ~r \in \mathbb{N}.$$
Most of the analysis of the numerical approximation in one dimension is easily extended to the two dimensional case as discussed by Morton et al. in \cite{Morton}. Truncation error of this approximation (\ref{3.2}) will remain as in  one dimensional case except some  additions due to the presence of $y$ variable \cite{Morton}. Usual analysis will give the CFL condition for stability in the following form
$$\frac{A \Delta t}{\Delta x}+\frac{B \Delta t}{\Delta x}\leq 1,$$ where $A$ and $B$ are the bounds of $a$ and $b$ respectively.\\
Proof of convergence follows  in similar way, leading to error in the approximation
 $$E^n\leq n\Delta t T_{max}\leq t_f T_{max},$$
 provided that the CFL condition is satisfied and $u$ has bounded derivatives of second order. \\
Similarly the  extension of the Leap-Frog numerical approximation (\ref{2.5}) with appropriate initial data and boundary-interval conditions is given by
\begin{equation}\label{3.3}
U_{i,j}^{n+1}=U_{i,j}^{n-1}-\frac{a_{i,j}^n\Delta t}{\Delta x}(U_{i+1,j}^n-U_{i-1,j}^{n})-\frac{b_{i,j}^n\Delta t}{\Delta y}(U_{i,j+1}^n-U_{i,j-1}^{n})+2\Delta t c_{i,j}^n U_{i-m_0,j-q_0}^n.
\end{equation}

\section{Illustration Examples}
Purpose of this section is to  include  some  test examples to validate the predicted  theory established in the paper and to illustrate the effect of point-wise delay  on the solution behavior. The maximum absolute errors for the Lax-Friedrichs approximation and  square root errors for the Leap-Frog approximation  are calculated using the double mesh principle~\cite{Doolan} as the exact solution for the considered examples are not available. We perform numerical computations using MATLAB.  The maximum absolute error is given by
$$E(\Delta x,\Delta t)=\max_{0\leq j\leq J,~0\leq n\leq N}{\left|U_{\Delta x}^{\Delta t}(j,n)-U_{\Delta x/2}^{\Delta t/2}(2j,2n)\right|}.$$
Similarly we find the square root errors for the Leap-Frog approximation.

\vskip 0.65cm \noindent
\(
\begin{array}{l}
\mathrm{Example\; 1.}~ \mathrm {We\; consider\;the\;differential\; difference\; equation\; (\ref{1.1})\; with\; the\; following\; coefficients\;}\\
\mathrm{and\; initial\mathrm{-}boundary\; conditions:} \\
\quad \: a(x,t)=\frac{1+x^2}{1+2xt+2x^2+x^4};
 \quad \quad \quad \quad\: b(x,t)=0.5;\\
\quad \: u(x,0)=\exp[-10(4x-1)^2];
\quad \quad \quad  \:   u(s,t)= 0,~\forall ~s\in[- \alpha ,0].
\end{array}
\)\\

\noindent We consider $\Omega=(0,1)$, $\Delta x=.001=\Delta t$. The numerical solution is plotted in Figure 1 with $\alpha=0.02$ at the time $t=0.5$ by both the numerical approximations. We observe that graphs of computed solution are very close to each other but there is a slight difference near the peak of the graphs. In Figure 2 and 3, we show the change in solution with the time by both the numerical approximations. As time increases, both the graphs shift to the right side with the time. The error tables 1 and 2  illustrate that the proposed numerical methods are convergent in both space and time direction. Error Tables are plotted by refining the grid points .

\vskip 0.65cm \noindent \(
\begin{array}{l}
\mathrm{Example\; 2.}~ \mathrm {We\; consider\;the\;differential\; difference\; equation\; (\ref{1.1})\; with\; the\; following\;  variable\; }\\
\mathrm{ coefficients\; and\; initial\mathrm{-}boundary\; conditions:} \\
\quad \: a(x,t)=\frac{1+x^2}{1+2xt+2x^2+x^4};
 \quad \quad \quad \quad\: b(x,t)=\frac{1}{1+x^2 t^2};\\
\quad \: u(x,0)=\exp[-10(4x-1)^2];
\quad \quad \quad  \:   u(s,t)= 0,~\forall ~s\in[- \alpha ,0].
\end{array}\)\\

\noindent In this example, we consider the original problem (\ref{1.1}) with variable  coefficients. We consider $\Omega=(0,1)$,  $\Delta x=.001=\Delta t$. The computed approximate solution is plotted in Figure 4 with $\alpha=0.05$ at time $t=0.5$ by both the numerical approximations. Both the approximations has slightly different at maxima but at other points graphs are very close to each other.   To show the effect of point-wise delay  on the solution behavior, we show the numerical solution with spatial variable $x$ in Figure 5 and 6 by the Lax-Friedrichs and the Leap-Frog approximation respectively. By changing the value of point-wise delay, we observe that as the value of $\alpha$ is increased, the height of impulse is decreased and width is decreased.  As the exact solution of this problem is not available, we calculate the maximum and square errors by refining the grid points.  Analysis of Tables 3 and 4 also verify the  convergence of both the approximations  in the space as well as time.

\vskip 0.65cm \noindent \(
\begin{array}{l}
\mathrm{Example\; 3.}~ \mathrm {We\; consider\;the\; 2-D\;differential\; difference\; equation\; (\ref{3.1})\; with\; the\; following\;}\\
\mathrm{ coefficients\; and\; initial\mathrm{-}boundary\; conditions:} \\
\quad \: a(x,y,t)=\frac{1+x^2+y^2}{1+2(x+y)t+2(x^2+y^2)+x^4};
 \quad \quad \quad \: b(x,y,t)=\frac{1}{1+(x^2+y^2) t^2};
  \quad \quad \quad \:c(x,y,t)=0.1;\\
\quad \: u(x,0)=\exp[-10(4x+4y-1)^2];\\
\quad \:   u(s_1,s_2,t)= 0,~\forall ~s_1\in[- \alpha ,0]\quad  \mathrm{and} \quad\forall ~s_2\in[- \beta ,0].
\end{array}\)\\

\noindent  We consider the two dimensional  problem (\ref{3.1}) with variable  coefficients. We consider $\Omega=(0,1)X(0,1)$,  $\Delta x= \Delta y=.01$ and time step $\Delta t=.001$. The  approximate numerical solutions are plotted  with $\alpha=0.5$ and $\beta=0.5$ at time $t=0.5$ by Lax-Friedrichs and Leap-Frog scheme in Figure 7 and Figure 8  respectively.
\clearpage

\begin{table}
\caption{The maximum absolute errors for Example 1 by using Lax-Friedrichs scheme with $\alpha=0.05$}
\begin{center} \footnotesize
\begin{tabular}{c c c c c c}\hline
  $\Delta t\downarrow$ $\Delta x \rightarrow$ &$1/100$   & $1/200$  & $1/400$  & $1/800$ \\ \hline
$\Delta x/2$ &$0.053623$	&$0.024289$	&$0.011064$ &$0.004520$ \\
$\Delta x/4$ &$0.026758$	&$0.011382$	&$0.006565$	&$0.002265$ \\
$\Delta x/8$  &$0.014642$	&$0.005163$	&$0.002344$	&$0.001041$   \\
$\Delta x/16$  &$0.008525$	&$0.004416$	&$0.001789$	&$0.000836$  \\ \hline
\end{tabular}
\end{center}
\end{table}
\begin{table}
\caption{The  square root errors for Example 1 by using Leap-Frog scheme with $\alpha=0.5$}
\begin{center} \footnotesize
\begin{tabular}{c c c c c c}\hline
  $\Delta t\downarrow$ $\Delta x \rightarrow$  &$1/100$   & $1/200$  & $1/400$  & $1/800$ \\ \hline
$\Delta x/2$ &$0.039607$	&$0.018968$	&$0.009380$	 &$0.003678$  \\
$\Delta x/4$ &$0.017185$	&$0.007408$	&$0.003681$	&$0.001538$   \\
$\Delta x/8$  &$0.008517$	&$0.003695$	&$0.001840$	&$0.000769$   \\
$\Delta x/16$  &$0.003750$	&$0.001347$	&$0.000770$	&$0.000285$   \\ \hline
\end{tabular}
\end{center}
\end{table}
\begin{table}
\caption{The maximum absolute errors for Example 2 by using Lax-Friedrichs scheme with $\alpha=0.05$}
\begin{center} \footnotesize
\begin{tabular}{c c c c c c}\hline
  $\Delta t\downarrow$ $\Delta x \rightarrow$  &$1/100$   & $1/200$  & $1/400$  & $1/800$ \\ \hline
$\Delta x/2$ &$0.053308$	&$0.024327$	&$0.011097$	&$0.006541$ \\
$\Delta x/4$ &$0.027900$	&$0.011334$	&$0.005577$	&$0.002174$  \\
$\Delta x/8$  &$0.017055$	&$0.004971$	&$0.002233$	&$0.001044$  \\
$\Delta x/16$  &$0.013255$	&$0.002274$	&$0.001042$	&$0.000433$  \\ \hline
\end{tabular}
\end{center}
\end{table}
\begin{table}
\caption{The  square root errors for Example 2 by using Leap-Frog scheme with $\alpha=0.1$}
\begin{center} \footnotesize
\begin{tabular}{c c c c c c}\hline
  $\Delta t\downarrow$ $\Delta x \rightarrow$ &$1/100$   & $1/200$  & $1/400$  & $1/800$ \\ \hline
$\Delta x/2$ &$0.041308$	&$0.019827$	&$0.009815$	&$0.004897$  \\
$\Delta x/4$ &$0.020042$	&$0.009843$	&$0.004901$	&$0.002448$  \\
$\Delta x/8$  &$0.009951$	&$0.004915$	&$0.002450$	&$0.001224$  \\
$\Delta x/16$  &$0.004969$	&$0.002457$	&$0.001225$	&$0.000612$  \\ \hline
\end{tabular}
\end{center}
\end{table}

\clearpage
\begin{figure}[h]
\includegraphics{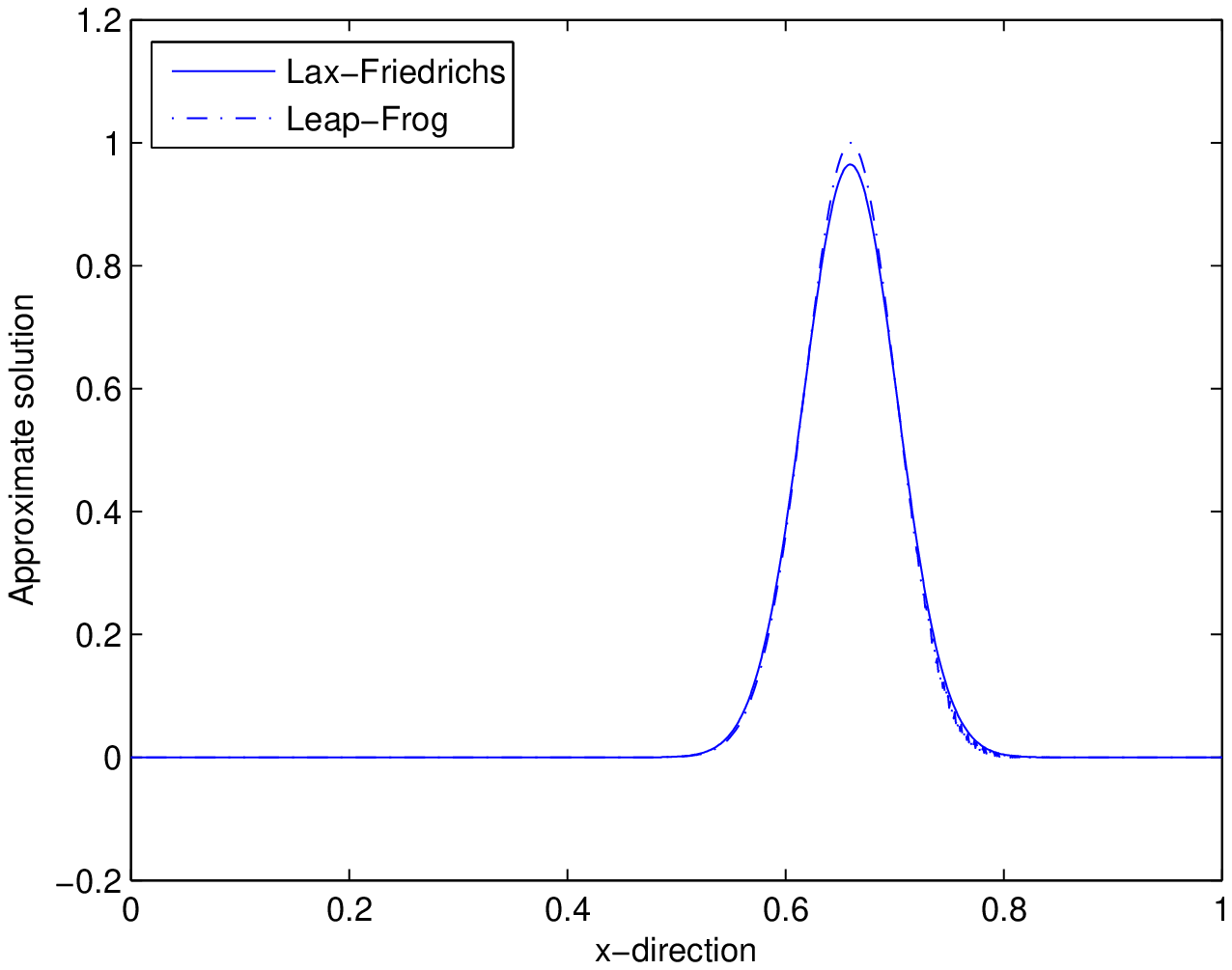}
\caption{Example 1: The approximate  solution  by both the schemes for $\alpha=0.8$ at $t=0.5$}
\includegraphics{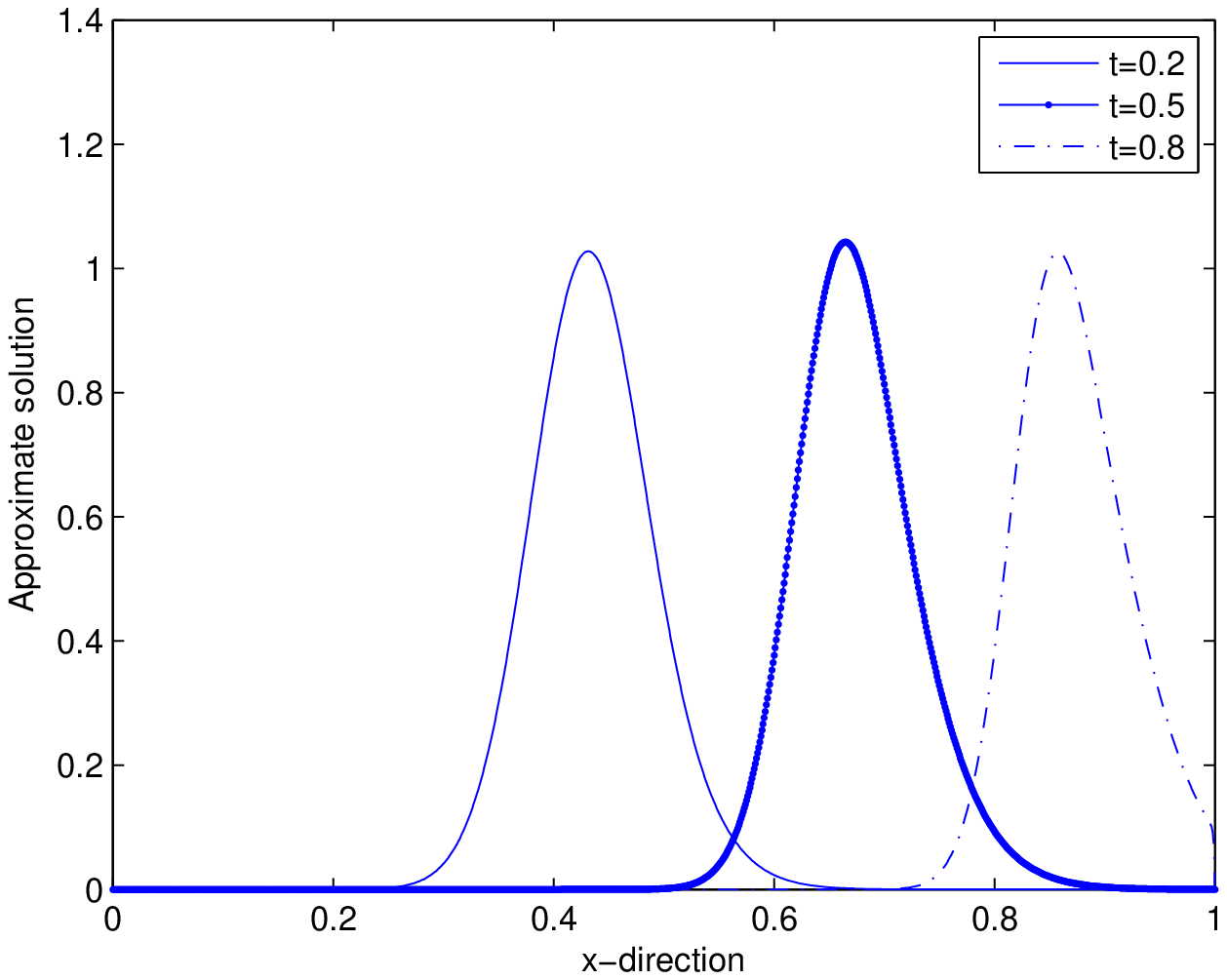}
\caption{Example 1: The approximate  solution at different time levels  with the Lax-Friedrichs scheme for  $\alpha=0.08$}
\end{figure}

\begin{figure}[h]
\includegraphics{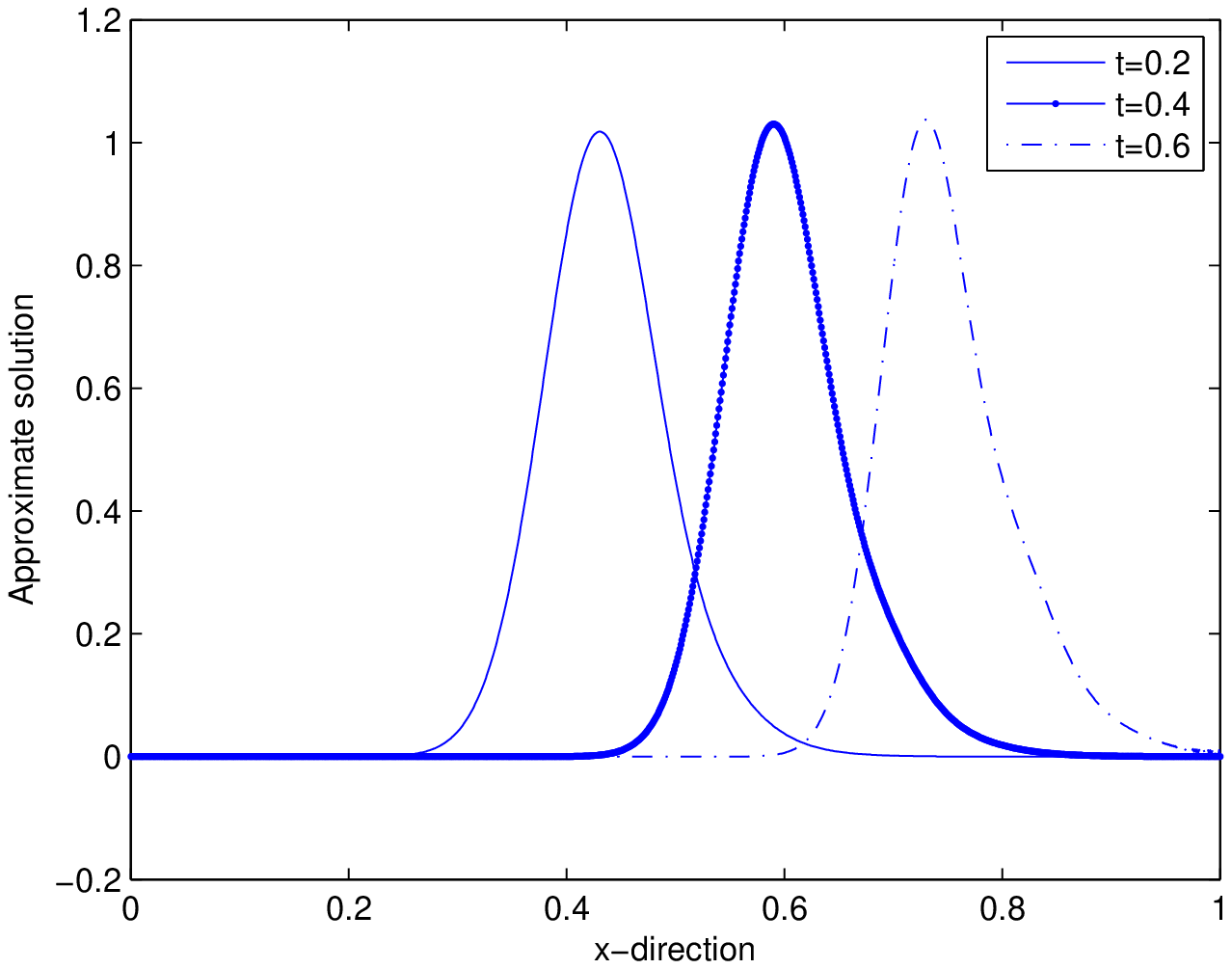}
\caption{Example 1: The approximate  solution at different time levels  with the Leap-Frog scheme for  $\alpha=0.1$}
\includegraphics{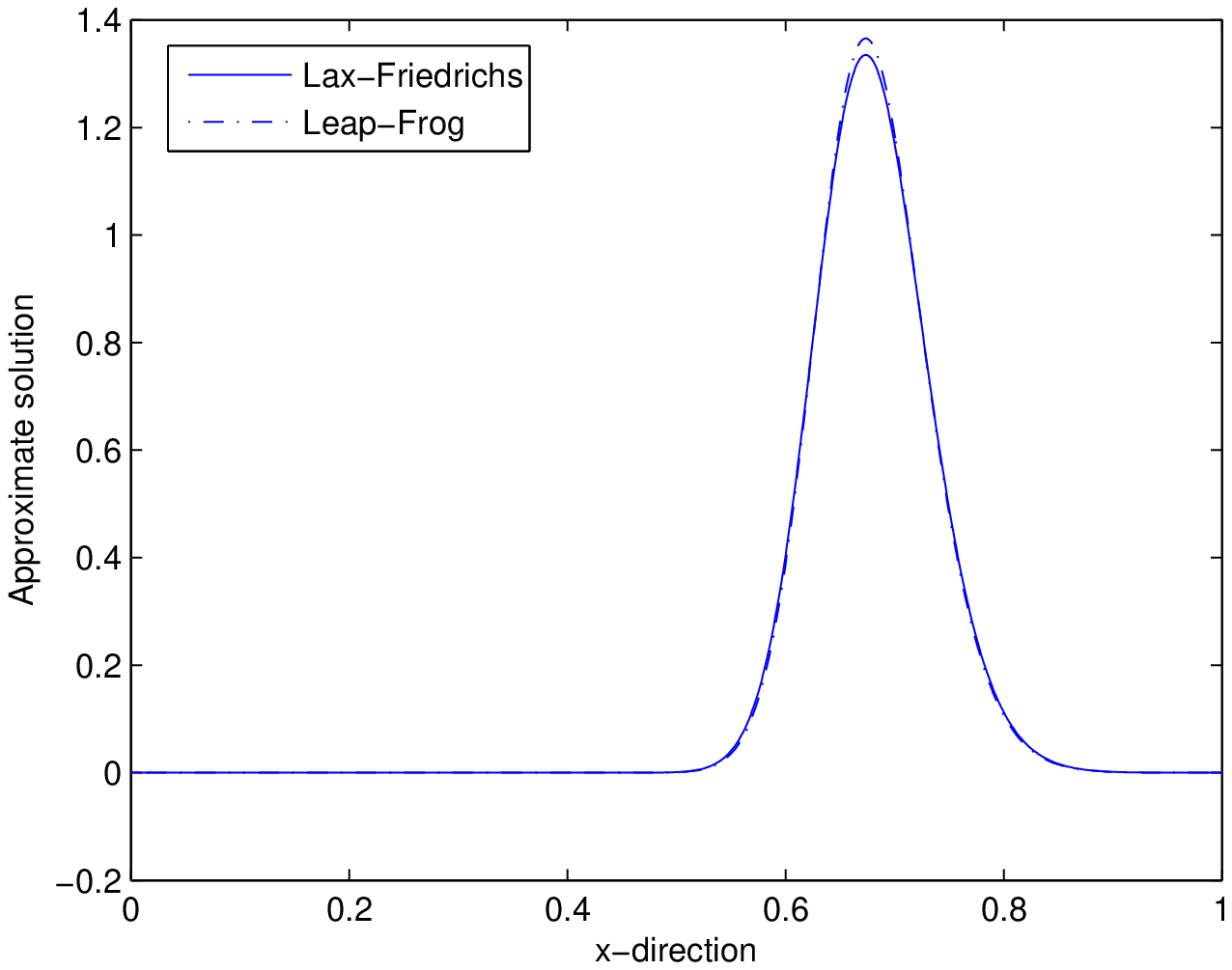}
\caption{Example 2: The approximate  solution  by both the schemes for $\alpha=.05$ at $t=0.5$}
\end{figure}

\begin{figure}[h]
\includegraphics{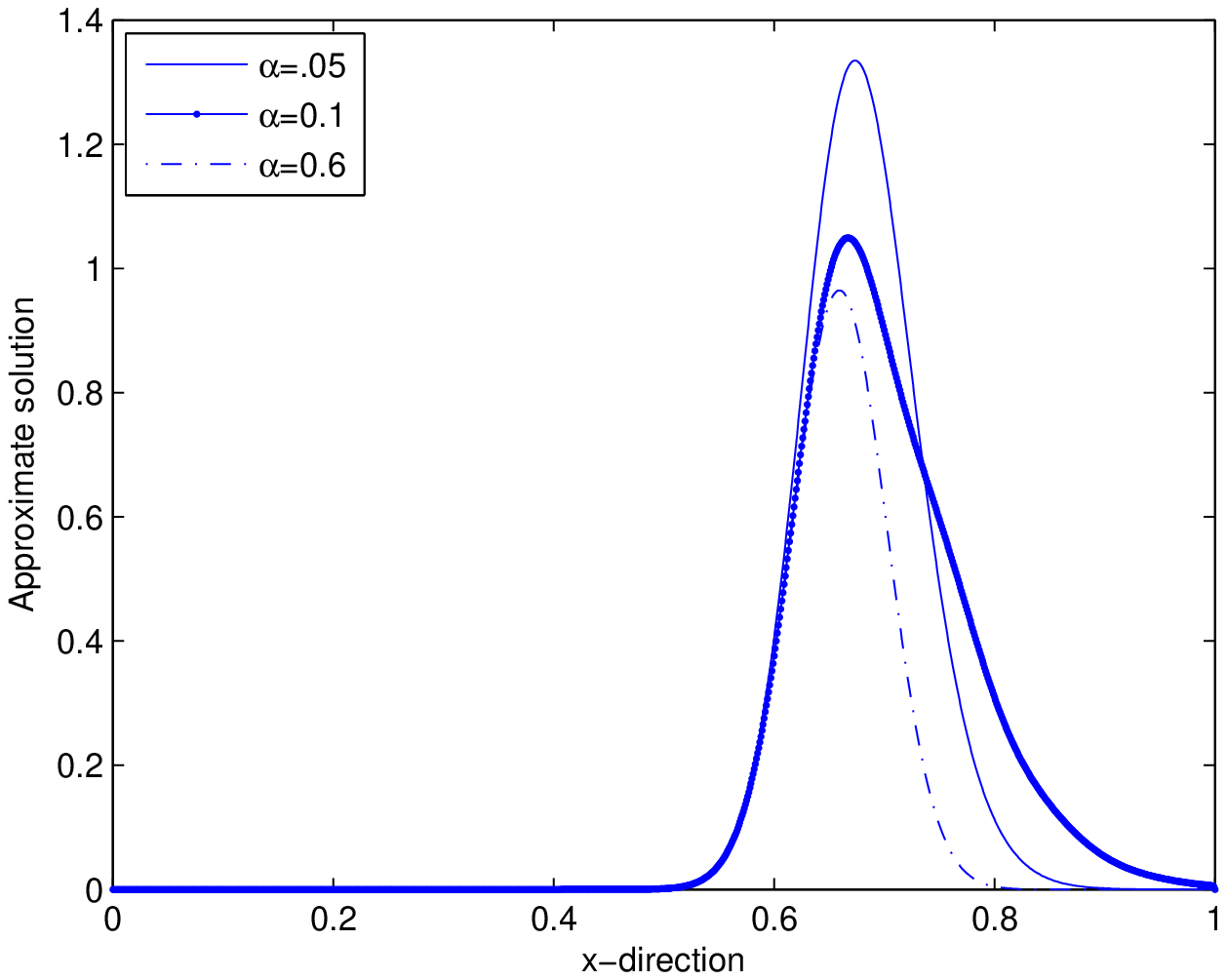}
\caption{Example 2: The effect of the point-wise delay on solution at $t=0.5$ by Lax-Friedrichs scheme}
\includegraphics{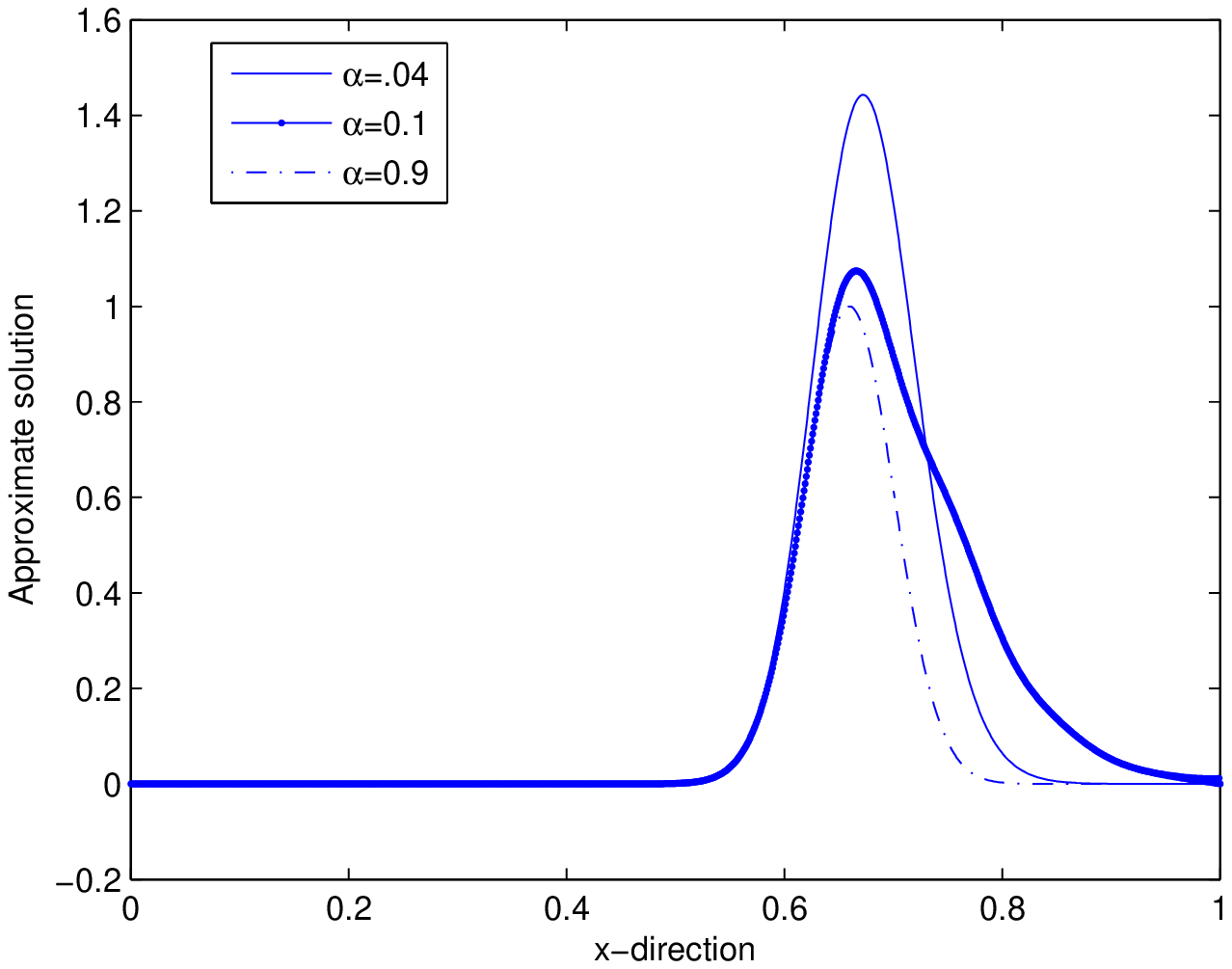}
\caption{Example 2: The effect of the point-wise delay on solution at $t=0.5$ by Leap-Frog scheme}
\end{figure}

\begin{figure}[h]
\includegraphics{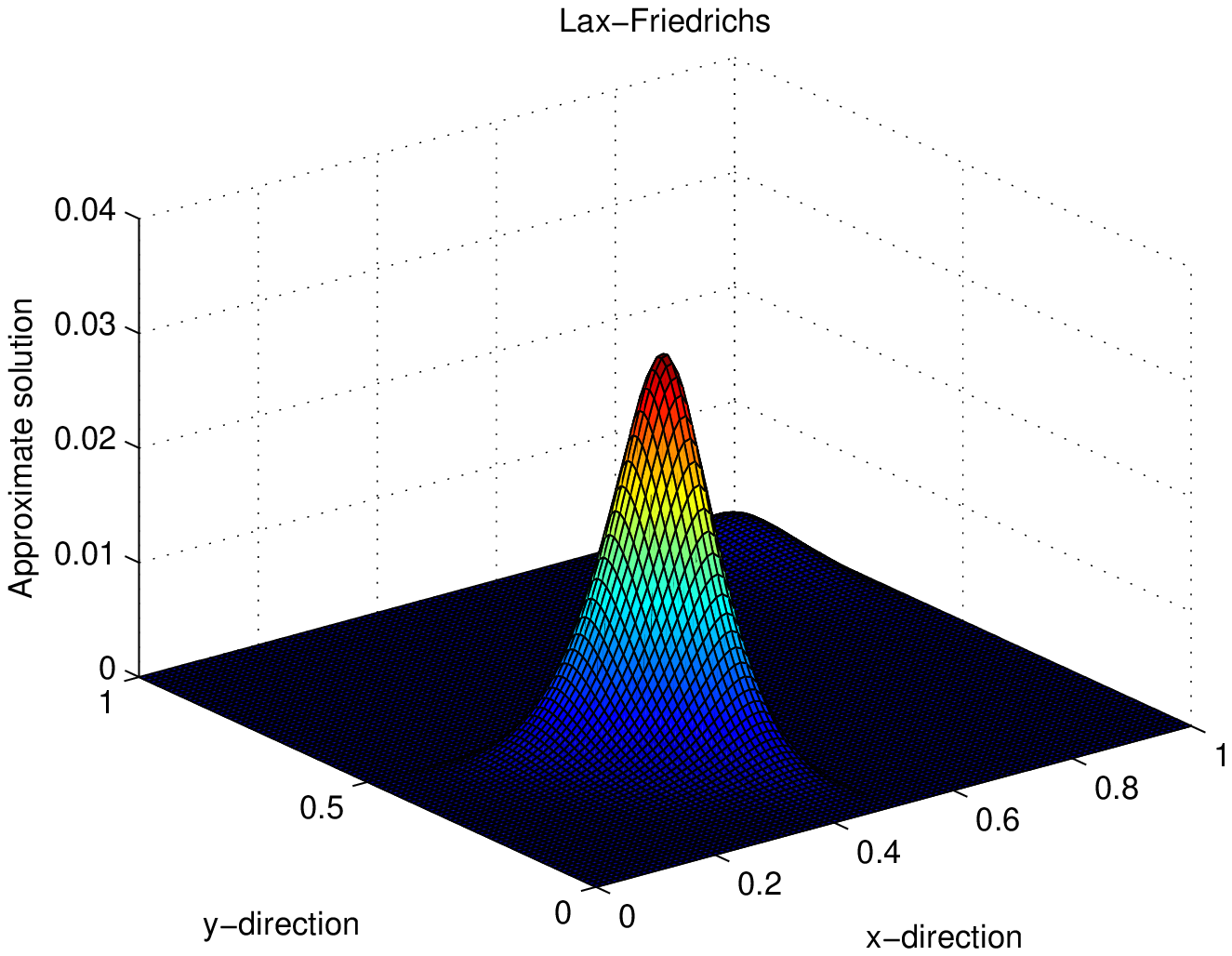}
\caption{Example 3: The approximate  solution  for $\alpha=0.5$ and $\beta=0.5$ at $t=0.5$}
\includegraphics{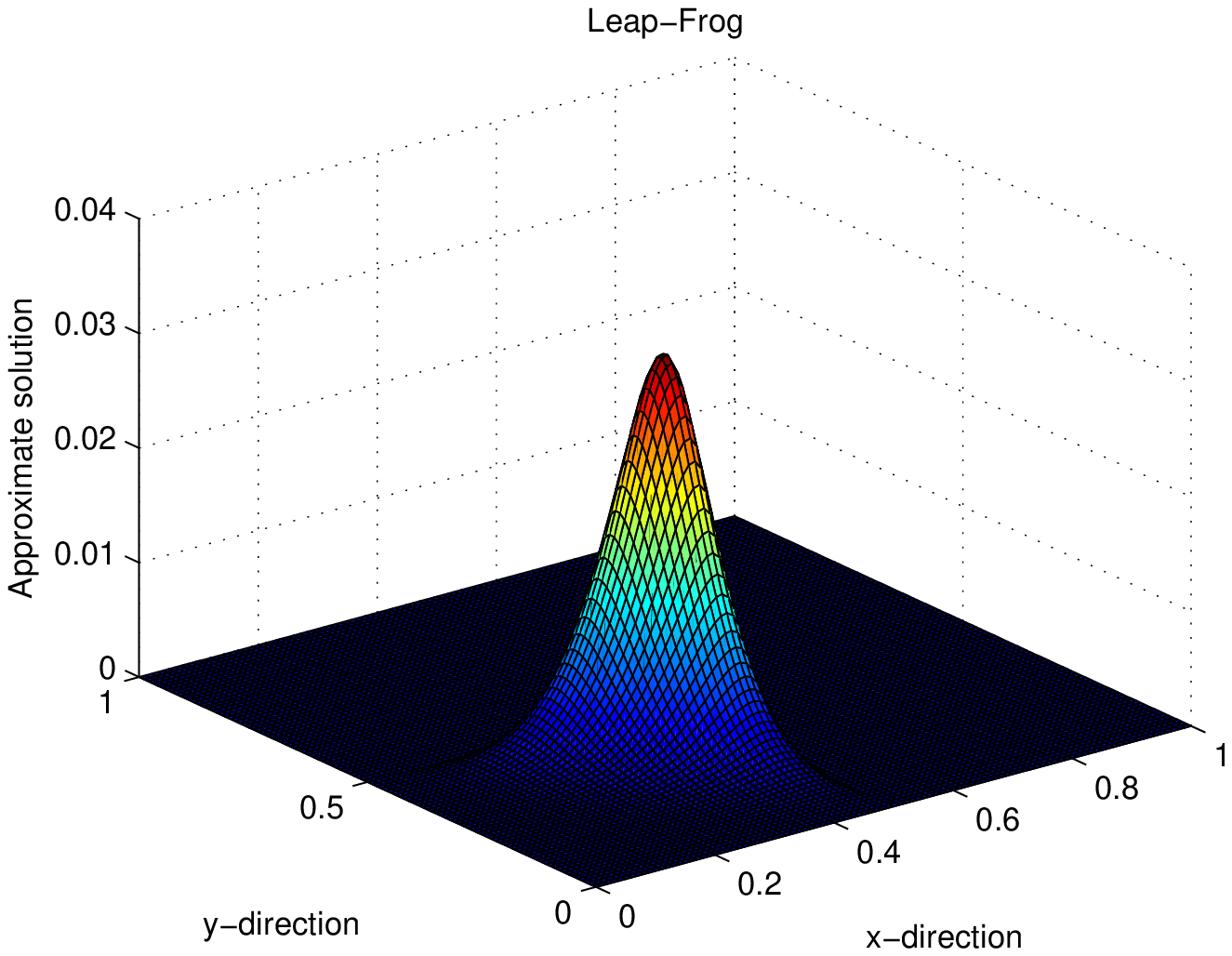}
\caption{Example 3: The approximate  solution  for $\alpha=0.5$ and $\beta=0.5$ at $t=0.5$}
\end{figure}

\clearpage
\section{Concluding Remarks}
We proposed in this paper a way to process  hyperbolic partial differential equation with point-wise delay.  It mainly relies on the constructions of finite difference approximations of order greater than one in space. In this numerical analysis, two explicit  numerical methods based on the Lax-Friedrichs and the Leap-Frog finite difference methods are constructed to find the numerical  solution of target problem with point-wise delay. The consistency, stability  and convergence analysis proves that proposed numerical schemes are consistent, stable with CFL condition and convergent in both space and time.  These second-order numerical methods in space  maintains the height and width better than a first-order scheme as discussed by the authors in \cite{Singh2}. The effect of point-wise delay on the solution behavior is shown by taking some test examples. Error Tables illustrate the  fact that methods are  convergent in space and time. Also we extends our ideas in higher space dimensions and include numerical experiments to show the behavior of solution in two space dimensions. \\
Finally, we remark that the strategy developed here can be applied to a problem having multiple point-wise delay or advance or both.\\

\noindent\textbf{Acknowledgments:} The first author is thankful to the Council of Scientific and Industrial Research, New Delhi, India for providing financial assistance in terms of Senior Research Fellowship. The authors are also thankful to Professor Bernardo Cockburn, School of Mathematics, University of Minnesota, Minneapolis, USA for his invaluable suggestions while preparing this manuscript.
%The anonymous reviewers are acknowledged with thanks for their helpful comments to improve the manuscript.

\end{document}